\documentclass[11pt]{amsart}

\usepackage{amsmath}
\usepackage{amssymb}
\usepackage{graphicx}
\newtheorem{theorem}{Theorem}[section]

\newtheorem{proposition}[theorem]{Proposition}
\theoremstyle{definition}

\numberwithin{equation}{section}

\title[Accretivity and form boundedness]{Accretivity and form boundedness  of \\second order differential operators}

\author[V. G. Maz'ya]{Vladimir G. Maz'ya}

\address{Department of Mathematics, Link\"oping University, SE-581 83, Link\"oping,
Sweden and RUDN University, 6 Miklukho-Maklay St., Moscow, 117198, Russia}
\email{{\tt vlmaz@mai.liu.se}}

\author[I. E. Verbitsky]{Igor E. Verbitsky}

\address{Department of Mathematics, University of Missouri, Columbia, Missouri 65211, USA}
\email{{\tt verbitskyi@missouri.edu}}

\keywords{Accretivity, form boundedness, general second order differential operators}

\subjclass[2010]{Primary 35J15,  42B37; Secondary 31B15, 35J10}

\newtheorem*{main1}{{\bf Theorem I}}
\newtheorem*{main2}{{\bf Theorem II}}
\newtheorem*{main3}{{\bf Theorem III}}
\newtheorem*{main4}{{\bf Theorem IV}}
\newtheorem*{main5}{{\bf Theorem V}}
\newtheorem*{main6}{{\bf Theorem VI}}

\newcommand{\R}{{\mathbb R}}

\thanks{The first author is supported in part  by the ``RUDN University Program 5-100"}

\dedicatory{In memory of Aizik Volpert}

\begin{document}

\begin{abstract} 
Let  $\mathcal{L}$ be the general second order differential operator   
with complex-valued 
distributional coefficients $A=(a_{jk})_{j, k=1}^n$, $\vec{b}=(b_{j})_{j=1}^n$, and $c$  in an open set $\Omega \subseteq \R^n$ ($n \ge 1$), with principal part  either   in the divergence form,
 $\mathcal{L} u=  {\rm div} \, (A \nabla u) +  \vec{b} \cdot\nabla u + c \, u$, or 
 non-divergence form, 
$
\mathcal L u= \sum_{j, \, k=1}^n \, a_{jk} \, \partial_j \partial_k u  +  \vec{b} \cdot\nabla u + c \, u
$.

We give a survey of the results by the authors  which characterize the following two properties 
of  $\mathcal{L}$:  

(1) $-\mathcal{L}$ is  
 accretive, i.e., ${\rm Re} \, \langle -\mathcal L u, 
 \, u\rangle \ge 0$; 
 
 (2) $\mathcal L$ is  form bounded, i.e., 
$\vert  \langle \mathcal L u, u \rangle  \vert \le C \, \Vert \nabla u \Vert_{L^2(\Omega)}^2$, 

 for all complex-valued $u \in C^\infty_0(\Omega)$.

\end{abstract}

\maketitle

\tableofcontents

\section{Introduction}\label{Introduction}

We consider   the general second order 
differential operator  
\begin{equation}\label{E:0.0}
\mathcal{L}_0 u = \sum_{j, \, k=1}^n \, a_{jk} \, \partial_j \partial_k u +  
\sum_{j=1}^n \, b_{j} \,  \partial_j u + c \, u,  
\end{equation}
in an open set $\Omega \subseteq \R^n$, with $A=(a_{jk}) \in D'(\Omega)^{n\times n}$,  $\vec{b} =(b_{j})\in D'(\Omega)^n$, and $c\in D'(\Omega)$,  where 
 $D'(\Omega)=C^\infty_0(\Omega)^*$  is the space of complex-valued distributions 
 in $\Omega$.  

We discuss the \textit{accretivity} 
property of  $-\mathcal{L}_0$ (or, equivalently,   \textit{dissipativity} of  
$\mathcal{L}_0$), i.e., 
 \begin{equation}\label{accr'}
{\rm Re} \,  \langle -\mathcal{L}_0 \, u, \, u \rangle \ge 0,  
\end{equation}
for all complex-valued functions  $u \in C^\infty_0(\Omega)$.

More general   differential operators 
\begin{equation}\label{E:0.b*}
\mathcal{L}_1 u=  \sum_{j, \, k=1}^n \, a_{jk} \, \partial_j \partial_k u +  \vec{b_1} \cdot\nabla u + {\rm div}\, 
(  \vec{b_2} \, u) + c_1 \, u, 
\end{equation}
with $\vec{b_1}, \vec{b_2} \in D'(\Omega)^n$ and $c_1\in D'(\Omega)$ 
can be treated as well, since $\mathcal{L}_1$ is immediately  reduced to 
$\mathcal{L}_0$ with $\vec{b} =\vec{b_1} + 
\vec{b_2}$ and $c=c_1 + {\rm div}\,  \vec{b_2}$.

Our main results on the accretivity problem for general differential operators are discussed in 
Sec. \ref{accretivity} below. (See  
Propositions \ref{prop_1} and  \ref{prop_2}, as well as  Theorem V for $n=1$, and Theorem VI for $n \ge 2$.)

For the sake of simplicity, we will focus in the Introduction on the operator 
$$\tilde{\mathcal{L}}=\Delta + \vec{b} \cdot \nabla +c,$$ 
whose principal part is the Laplacian $\Delta$, and the coefficients 
  $\vec{b}=(b_j)$ and $c $ are locally integrable functions in $\R^n$. Then 
  the sesquilinear form of $-\tilde{\mathcal{L}}$ is given by 
  \begin{align}\label{sesq-l}
\langle - \tilde{\mathcal{L}} u, v\rangle & =
 \int_{\R^n}  (  \nabla u \cdot \nabla \overline{v}  - \vec{b} \cdot\nabla u \, \, 
\overline{v} - c \, u \,  \overline{v} ) \, dx,      
\end{align}
where $u, \, v  \in C^\infty_0(\R^n)$.

In this special case, let 
\begin{equation}\label{sigma-1}
q ={\rm Re} \, c - \tfrac {1}{2}{\rm div} \, ({\rm Re} \, \vec{b}), \qquad \vec{d} =  
\tfrac{1}{2} ( {\rm Im} \,\vec{b}). 
\end{equation}

We denote by $\mathcal{H}=\Delta + q$ the corresponding Schr\"{o}dinger operator. The quadratic form associated with 
$-\mathcal{H}$ in the case $q \in L^1_{{\rm loc}}(\R^n)$ is given by 
\begin{equation}\label{norm-1}
[h]_{\mathcal{H}}^2 := \langle  -\mathcal{H} h, h\rangle= 
\int_{\R^n}  ( |\nabla h|^2  - q \,  |h|^2 ) \, dx, \quad  h \in C^\infty_0(\R^n). 
\end{equation}

\begin{main1}\label{Theorem I}   {\it Let $\tilde{\mathcal{L}}=\Delta + \vec{b} \cdot \nabla +c$, where ${\rm Re} \, \vec{b}\in W^{1, 1}_{{\rm loc}}(\R^n)$, 
and ${\rm Im} \, \vec{b}, \, c \in L^1_{{\rm loc}}(\R^n)$. 
Let $q$, $\vec{d}$ be  given by 
 \eqref{sigma-1}. 
 Then 
the operator $-\tilde{\mathcal{L}}$ is accretive 
if and only if  the following two conditions hold: 

   {\rm (i) } The operator $-\mathcal{H}$ is nonnegative definite, i.e., 
\begin{equation}\label{norm-2}
[h]_{\mathcal{H}}^2=\int_{\R^n}  ( |\nabla h|^2  -  q \, |h|^2 ) \, dx \ge 0, \end{equation}
for all real (or complex-valued)  $h\in C^\infty_0(\R^n)$.

   {\rm (ii) }  The commutator 
inequality 
\begin{equation}\label{comm-1}
\left \vert \int_{\R^n} \vec{d} \cdot(u \nabla v -v \nabla u) \, dx  \right\vert 
\le  [u]_{\mathcal{H}} \, [v]_{\mathcal{H}}   
\end{equation}
holds for all real-valued  $u, v \in C^\infty_0(\R^n)$.}
\end{main1}
\smallskip

A necessary and sufficient condition for property \eqref{norm-2} was obtained 
in \cite[Proposition 5.1]{JMV1} (see Sec. \ref{sec2.3} below). Concerning condition \eqref{comm-1}, we observe 
that, under the upper and lower bounds  on the quadratic form \eqref{norm-2} discussed 
in Sec. \ref{sec2.5}, the expressions  $[u]_{\mathcal{H}}$  and $[v]_{\mathcal{H}}$  on the right-hand side of  \eqref{comm-i}
 can be replaced,  up to a constant multiple, with the corresponding Dirichlet norms $|| \nabla u||_{L^{2} (\R^n)}$ and $|| \nabla v||_{L^{2} (\R^n)}$, respectively. 
 Then the corresponding commutator inequality
 \begin{equation}\label{comm-2}
\left \vert \int_{\R^n} \vec{d} \cdot(u \nabla v -v \nabla u) \, dx  \right\vert 
\le  C \, \, || \nabla u||_{L^{2} (\R^n)} \, 
|| \nabla v||_{L^{2} (\R^n)},  
\end{equation}
for all (real-valued or complex-valued) $u, v \in C^\infty_0(\R^n)$,  
can be characterized completely as follows (see \cite[Lemma 4.8]{MV3}).

\begin{main2}\label{Theorem II}{\it 
Let $\vec{d} \in L^1_{{\rm loc}}(\R^n)$, $n \ge 2$. Then inequality \eqref{comm-2}
 holds if and only if 
\begin{equation}\label{E:decomp}
\vec d = \vec c + {\rm Div} \, F,
\end{equation}
where  $F \in 
{\rm BMO}(\R^n)^{n \times n}$ is a skew-symmetric 
matrix field, and $\vec c$ satisfies the condition 
\begin{equation}\label{E:1.5ub}
\int_{\R^n} |\vec c|^2 \, |u|^2 \, dx  \leq \, C \,
 || \nabla u||^2_{L^{2} (\R^n)},   
\end{equation}
where the constant $C$ does not depend on $u \in C^\infty_0(\R^n)$.

Moreover, if \eqref{comm-2} holds, then  
 \eqref{E:decomp} is valid 
with $\vec c = \nabla \Delta^{-1} ({\rm div} \, \vec d)$ satisfying \eqref{E:1.5ub}, and 
$F = \Delta^{-1} ({\rm Curl} \, \vec d) \in {\rm BMO}(\R^n)^{n \times n}$. 

In the case $n=2$, necessarily $\vec c =0$, and  $\vec{d}= (-\partial_2 f, \partial_1 f)$
with $f \in {\rm BMO}(\R^2)$ in the above statements.} 
\end{main2} 
\smallskip

Here the gradient $\nabla$, and the matrix  operators ${\rm Div}$, ${\rm Curl}$  are  understood in the sense 
of distributions (see Sec. \ref{preliminaries}). Expressions $\Delta^{-1} ({\rm div} \, \vec d)$, $\Delta^{-1} ({\rm Curl} \, \vec d)$, etc., are defined in terms of the weak-$*$ BMO convergence (details can be found in \cite{MV3}, \cite{MV6}). Theorems I \&  II yield an explicit criterion of accretivity for $-\tilde{\mathcal{L}}$ (see Theorem VI below in the general case).

More general commutator inequalities related to compensated compactness theory \cite{CLMS}  
 were studied earlier by the authors \cite{MV3}
 in the framework of the 
\textit{form boundedness} problem,  
\begin{equation}\label{E:1.2old}
| \langle \mathcal{L}_0 \, u, \, v \rangle | \leq \, C \, || \nabla u||_{L^{2} (\R^n)} \, 
|| \nabla v||_{L^{2} (\R^n)},   
\end{equation}
where the constant $C$ does not depend on $u, \, v \in C^\infty_0(\R^n)$.

If \eqref{E:1.2old} holds, 
then 
$\langle \mathcal{L}_0 \, u, \, v \rangle$ can be extended by continuity to 
$u, v \in L^{1, \, 2} (\R^n)$ ($n \ge 3$). Here $L^{1, \, 2}(\R^n)$   is 
the completion of (complex-valued) $C^\infty_0(\R^n)$ functions with respect to the  norm 
$|| u||_{L^{1, \, 2}(\R^n)} = ||\nabla  u||_{L^2(\R^n)}$. 
 Equivalently,  
\begin{equation}\label{E:dual-op}
\mathcal{L}_0\!:  L^{1, \, 2}(\R^n)\to L^{-1, \, 2}(\R^n)
\end{equation}
 is a bounded operator, where  
$L^{-1, \, 2}(\R^n)=L^{1, \, 2}(\R^n)^*$ is a dual 
Sobolev space. Analogous problems have been studied in \cite{MV2}--\cite{MV5}  for the inhomogeneous Sobolev space 
 $W^{1, \, 2}(\R^n)$, fractional Sobolev spaces, 
 infinitesimal 
 form boundedness, and other related questions 
 (see Sec. \ref{form-boundedness} below).  
 
 In the special case of the operator $\tilde{\mathcal{L}}$, we have the following 
 characterization of form boundedness. 
 
 \begin{main3}\label{Theorem III} {\it Let 
$\tilde{\mathcal{L}} = \Delta +
\vec b \cdot \nabla + q$, where  
 $\vec b \in L^1_{{\rm loc}}(\R^n)^n$ and $q \in L^1_{{\rm loc}}(\R^n)$, $n \ge 2$.  
 Then 
the following statements hold. 

  {\rm (i) } The sesquilinear form of $\tilde{\mathcal{L}}$ given by \eqref{sesq-l} 
is bounded 
if and only  if  
 $\vec b$ and $q$ can be represented respectively in the form 
\begin{equation}\label{E:1.4ua}
\vec b = \vec c + {\rm Div} \, F, \qquad q = {\rm div} \, \vec h,
\end{equation} 
where $F$ 
 is a skew-symmetric matrix field such that
\begin{equation}\label{E:1.4va}
F \in {\rm BMO(\R^n)}^{n\times n},
\end{equation} 
whereas $\vec c$ and $\vec h$ satisfy the condition 
\begin{equation}\label{E:1.5ua}
\int_{\R^n} (|\vec c|^2 + | \vec h |^2) \, |u|^2 \, dx  \leq \, C \,
 || \nabla u||^2_{L^{2} (\R^n)},   
\end{equation}
where the constant $C$ does not depend on $u \in C^\infty_0(\R^n)$.

  {\rm (ii) } If  the sesquilinear form of $\tilde{\mathcal{L}}$ is bounded, then $\vec c$,  
$F$, and $\vec h$  in decomposition {\rm (\ref{E:1.4ua})} can be determined  
explicitly by    
\begin{align}\label{E:1.6ua}
\vec c & = \nabla  \Delta^{-1} ({\rm div} \, \vec b), 
 \qquad \vec h = 
 \nabla (\Delta^{-1} \, q), \\
\qquad F & = \Delta^{-1} ({\rm Curl} \, \vec b),  
\label{E:1.6uua}
\end{align} 
so that conditions \eqref{E:1.4va}, 
\eqref{E:1.5ua} hold.

If $n=2$, then \eqref{E:1.5ua} yields that  $\vec c=0$ and $\vec h=0$, so that $q=0$ and $\vec b=(-\partial_2 f, \partial_1 f)$
with $f \in {\rm BMO}(\R^2)$.}  
\end{main3}
\smallskip

The form boundedness problem \eqref{E:1.2old} for the 
general second order differential operator $\mathcal{L}_0$ in the case $\Omega =\R^n$ 
was  characterized by the authors in \cite{MV3} using harmonic 
analysis  and potential theory methods. 
These results are discussed in Sec. \ref{form-boundedness} below. We observe that no ellipticity assumptions 
are imposed  
on the principal part $\mathcal{A}$ of $\mathcal{L}_0$ in this context. 

For the Schr\"{o}dinger operator 
$\mathcal{H}=\Delta + q$ with $q \in D'(\Omega)$, where 
either $\Omega =\R^n$, or $\Omega$ is a bounded domain that supports 
Hardy's inequality (see \cite{An}), a characterization  of form boundedness was obtained earlier in \cite{MV1}. 
A different approach for $\mathcal{H}=  {\rm div} \, (P \nabla \cdot) + q$ in 
general open sets $\Omega\subseteq \R^n$,  under the uniform ellipticity assumptions on $P$, was developed 
in \cite{JMV1}. (We remark that these assumptions on $P$ can be relaxed in a substantial way.)  
There is also a quasilinear version for operators of the $p$-Laplace type 
(see \cite{JMV2}).

Both the accretivity and form boundedness properties have numerous applications.  They include  problems in mathematical quantum mechanics (\cite{RS}, \cite{RSS}),   PDE theory (\cite{EE}, \cite{GGM}, \cite{K}, 
\cite{KP}, \cite{M}, \cite{NU}, \cite{H}, \cite{Ph}), fluid mechanics and Navier-Stokes equations (\cite{GP}, \cite{KT}, \cite{SSSZ}, \cite{T}), semigroups  
and Markov processes  
(\cite{LPS}), homogenization theory (\cite{ZP}), harmonic analysis 
(\cite{CLMS}, \cite{FNV}), etc. 

We conclude the Introduction with the observation that, for 
the form boundedness property, the case of  complex-valued coefficients  is easily reduced to the real-valued case. 
 In contrast, for the accretivity property,  complex-valued coefficients  
lead to additional difficulties that appear when the matrix ${\rm Im} \, A$ is not symmetric, or the imaginary part of $\vec{b}$ is nontrivial. 

\section{Preliminaries}\label{preliminaries}

Let   $\Omega \subseteq \R^n$ ($n \ge 1$) be an open set.  
The matrix row divergence 
operator  
${\rm Div}\!\!:  D'(\Omega)^{n \times n}\to D'(\Omega)^{n}$ is defined 
on matrix fields $F = (f_{j k})_{j, k =1}^n \in D'(\Omega)^{n \times n}$ by ${\rm Div} \, F = \left ( \sum_{k=1}^n \, \partial_k \, 
f_{j k}  \right)_{j=1}^n \in D'(\Omega)^n$. If $F$ is 
skew-symmetric, i.e., $f_{jk} = -f_{kj}$, 
then we obviously have ${\rm div} \, ({\rm Div} \, F) = 0$.

The matrix curl operator ${\rm Curl}\!\!:  D'(\Omega)^n \to D'(\Omega)^{n \times n}$ 
is defined on vector fields $\vec{f}=(f_k)_{k=1}^n$ 
 by ${\rm Curl} \, \vec{f} = (\partial_j f_k-\partial_k f_j)_{j, k=1}^n$. Clearly, 
 ${\rm Curl} \, \vec{f}$ is always a skew-symmetric matrix field.

 It will be convenient to use the notion of {\it admissible measures}  
 $\mathfrak{M}^{1, \, 2}_+(\Omega)$, i.e., nonnegative locally finite Borel measures $\mu$ in $\Omega$ 
 which obey  the trace inequality    
\begin{equation}\label{E:1.tr}
\Big(\int_{\Omega}  
  |u|^2  \, d \mu\Big)^{\frac{1}{2}}  \le C \, 
|| \nabla u||_{L^{2}(\Omega)}, \qquad \textrm{for all} \, \, u \in C^\infty_0(\Omega), 
\end{equation} 
where the constant $C$ does not depend on $u$. The least embedding constant $C$ in \eqref{E:1.tr} 
will be denoted by $||\mu||_{\mathfrak{M}^{1, \, 2}_+(\Omega)}$.  For admissible measures $q (x) \, dx$ 
with nonnegative  density $q\in L^1_{{\rm loc}}(\Omega)$, we write $q \in \mathfrak{M}^{1, \, 2}_+(\Omega)$.

Several characterizations of $\mathfrak{M}^{1, \, 2}_+(\Omega)$ are known. They can be formulated    
in terms of capacities \cite{M} or Green energies \cite{FNV}, \cite{QV}, and, in the case $\Omega=\R^n$, in terms of 
local maximal estimates \cite{KS}, pointwise 
potential inequalities \cite{MV1}, or dyadic Carleson measures 
\cite{V} (see also \cite{MV3}, \cite{MV6}).

Suppose that  the principal part $\mathcal{A} u$ of the general differential operator is given 
in the divergence form, 
\begin{equation}\label{E:div-form}
\mathcal{A} u =  {\rm div} \, (A \nabla u), \quad u \in C^\infty_0(\Omega). 
\end{equation}
Then we consider the operator 
\begin{equation}\label{E:0.b}
\mathcal{L} u=  {\rm div} \, (A \nabla u) +  \vec{b} \cdot\nabla u + c \, u, 
\end{equation}
with distributional coefficients $A=(a_{jk})$, $\vec{b} =(b_j)$, and $c$. 
The corresponding sesquilinear form $\langle \mathcal{L} u, v\rangle$ is given by 
\begin{equation}\label{E:div-sesq}
\langle \mathcal{L} u, v\rangle = - \langle A \nabla u, \nabla v \rangle 
+ \langle \vec{b} \cdot\nabla u, v\rangle + \langle c \,  u, v\rangle, 
\end{equation}
where $u, v \in C^\infty_0(\Omega)$ are complex-valued.

We observe that if $\mathcal{L}_0$ is given in the non-divergence form \eqref{E:0.0}, then 
\[
\mathcal{L}_0 = \mathcal{L} - {\rm Div} \, {A} \cdot  \nabla.
\] 
 (See, for instance, 
\cite{KP}, \cite{MV6}.)  
 Hence, we 
 can  express $\langle \mathcal{L}_0 u, v\rangle$ in the form 
\eqref{E:div-sesq}, with  $ \vec{b}- {\rm Div} \, {A}$ 
in place of $ \vec{b}$, for 
distributional coefficients $A$ and  $\vec{b}$.

This means that, without loss of generality, we may treat the accretivity property   
\begin{equation}\label{accr}
{\rm Re} \,  \langle -\mathcal{L} \, u, \, u \rangle \ge 0,  \quad \textrm{for all} \, \, u \in C^\infty_0(\Omega),
\end{equation}
 for  the divergence form operator $\mathcal{L}$ given by \eqref{E:0.b}.

This problem is of substantial interest even in the real-variable case, where 
the goal is to characterize  operators $-\mathcal{L}$ with real-valued coefficients whose quadratic form is nonnegative definite, 
\begin{equation}\label{pos-def}
\langle -\mathcal{L} h, h\rangle \ge 0, \quad \textrm{for all real-valued} \, \, h \in C^\infty_0(\Omega).
 \end{equation} 
 In this case the operator $-\mathcal{L}$ is called \textit{nonnegative definite}. 

 In the special case of Schr\"{o}dinger operators  
 \begin{equation}\label{schro}
 \mathcal{H} u = {\rm div} \, (P \nabla u) + \sigma \, u,
  \end{equation}
 with real-valued  $P\in D'(\Omega)^{n \times n}$ and $\sigma \in D'(\Omega)$, a characterization of this property 
 was obtained earlier in \cite[Proposition 5.1]{JMV1} under the assumption that $P$ is uniformly elliptic, i.e., 
 \begin{equation}\label{ell}
m \, ||\xi||^2 \le  P(x) \xi \cdot \xi \le M \, ||\xi||^2, \quad \textrm{for all} \, \, \xi \in \R^n, 
\, \,  \textrm{a.e.} \, \, x \in \Omega, 
 \end{equation}
 with the ellipticity  constants $m>0$ and $M<\infty$. 
  
An analogous characterization of   \eqref{pos-def} 
 for more general operators which include drift terms, $\mathcal L 
= {\rm div}  ( P \nabla \cdot) +
\vec{b} \cdot \nabla + c$, with real-valued coefficients and $P$ satisfying \eqref{ell},  
is given in Proposition \ref{prop_2}  below.

For the general differential operator 
in the form 
\eqref{E:div-form}, 
 we define the symmetric part $A^s$, and  co-symmetric (or skew-symmetric) part $A^c$,  
respectively, by 
 \begin{equation}\label{trans}
A^s =\frac{1}{2}(A + A^{\perp}), \quad A^c =\frac{1}{2}(A - A^{\perp}). 
 \end{equation}
 Here $A=(a_{jk})\in D'(\Omega)^{n \times n}$, and $A^{\perp}=(a_{kj})$ is the transposed matrix. 

For $-\mathcal{L}$ to be accretive, the matrix 
 $A^s$  must have a nonnegative definite real part: 
$P = {\rm Re} \, A^s$  should satisfy 
\begin{equation}\label{matr-def}
P \xi\cdot \xi \ge 0 \quad \textrm{for all} \, \, \xi \in \R^n, \quad \textrm{in} \, \, D'(\Omega). 
\end{equation}
Moreover, if the corresponding Schr\"{o}dinger operator $\mathcal{H}$ is defined by 
\eqref{schro} 
with 
\[
P= {\rm Re} \, A^s, \quad \sigma ={\rm Re} \, c - \frac {1}{2}{\rm div} \, ({\rm Re} \, \vec{b}), 
\]
then   $-\mathcal{H}$ must be 
nonnegative definite: 
\begin{equation}\label{norm}
[h]_{\mathcal{H}}^2 =\langle - \mathcal{H} h, h\rangle= 
\langle P \nabla h, \nabla h\rangle - \langle \sigma h, h\rangle \ge 0, 
\end{equation}
for all real-valued (or complex-valued) $h \in C^\infty_0(\Omega)$.

The rest of the accretivity problem for $\mathcal L$ (see Sec. \ref{sec2.1}) 
is reduced to the commutator 
inequality 
\begin{equation}\label{comm-i}
\left \vert \langle \vec{d},  u \nabla v -v \nabla u\rangle \right\vert 
\le  [u]_{\mathcal{H}} \, [v]_{\mathcal{H}},  
\end{equation}
for all real-valued  $u, v \in C^\infty_0(\Omega)$, where the real-valued vector field 
$\vec{d}$ is given by 
\begin{equation}\label{comm-d}
\vec{d}=  \tfrac{1}{2} \,  [ {\rm Im} \,\vec{b} - {\rm Div} ({\rm Im} \, A^c)].
\end{equation}

As mentioned in the Introduction, under some mild restrictions on $\mathcal{H}$, the ``norms''  $[u]_{\mathcal{H}}$  and $[v]_{\mathcal{H}}$  on the right-hand side of  \eqref{comm-i}
 can be replaced,  up to a constant multiple, with the corresponding Dirichlet norms $|| \nabla \cdot||_{L^{2} (\Omega)}$. This leads to explicit criteria of accretivity, such   as  Theorem VI below in the case $\Omega=\R^n$.

\section{Form boundedness}\label{form-boundedness}

We start with a  discussion of form boundedness for the general second order 
differential operator $\mathcal L$ in the form \eqref{E:0.b}, where $a_{ij}$,  $b_i$, and $c$ are real- or 
complex-valued distributions, on the homogeneous 
Sobolev space $L^{1, \, 2}(\R^n)$, and its inhomogeneous counterpart 
$W^{1, \, 2}(\R^n)$, obtained in \cite{MV3}. 

In particular, this leads to criteria of the \textit{relative form 
boundedness}   of the operator $\vec b \cdot \nabla + q$ with distributional coefficients $\vec b$ and $q$  
with respect to the Laplacian  $\Delta$ on $L^2(\R^n)$. 
Invoking  the so-called 
KLMN Theorem (see \cite[Theorem IV.4.2]{EE}; \cite[Theorem X.17]{RS}), we   
can then demonstrate 
that $\tilde{\mathcal{L}} = \Delta + \vec b \cdot \nabla + q$ is well defined,  
under appropriate 
smallness assumptions on 
$\vec b$ and $q$, as an 
m-sectorial operator on $L^2(\R^n)$. In this case, the  quadratic form domain 
of $\tilde{\mathcal{L}}$ 
coincides with $W^{1, \, 2}(\R^n)$.  
 
This yields a  characterization of the relative 
form boundedness  for the magnetic Schr\"odinger operator 
\begin{equation}\label{E:1.1m}
\mathcal{M} = (i \, \nabla + \vec a )^2 + q, 
\end{equation}
with  arbitrary vector potential $\vec a \in L^2_{{\rm loc}} (\R^n)^n$, and 
$q \in D'(\R^n)$ on $L^2(\R^n)$ with respect to $\Delta$ (see \cite{MV3}).

Our  approach is based on factorization of functions in Sobolev spaces and  integral estimates of 
potentials of equilibrium measures, combined with compensated compactness 
arguments, commutator estimates, and the idea of gauge invariance. Moreover,  
 an explicit Hodge decomposition is established for form bounded 
vector fields in $\R^n$. In this decomposition,  the irrotational part of the vector field 
is subject to a stronger restriction than its divergence-free counterpart. 

\subsection{Form boundedness in the homogeneous Sobolev space}\label{homogeneous} 

As was mentioned above, without loss of generality we may assume  
that the principal 
part of the differential operator is  in 
the divergence form, i.e.,  $\mathcal L= {\rm div} \, (A \, \nabla\cdot) +
\vec b \cdot \nabla + q$.

We present necessary and sufficient conditions on $A$, $\vec b$, 
and $q$, obtained in \cite[Theorem I]{MV3}, 
which ensure the boundedness in the homogeneous Sobolev space 
$L^{1, \, 2} (\R^n)$ of  
 the sesquilinear form associated with 
 $\mathcal L$:
\begin{equation}\label{E:1.1a}
| \langle \mathcal{L} \, u, \, v \rangle | \leq \, C \, ||u||_{L^{1, \, 2} (\R^n)} \, 
||v||_{L^{1, \, 2} (\R^n)},    
\end{equation}
where  $C$ does not depend on $u, \, v  \in C^\infty_0(\R^n)$, and 
$|| u||_{L^{1, \, 2}(\R^n)} = ||\nabla  u||_{L^2(\R^n)}$.

\begin{main4}\label{Theorem IV} {\it Let 
$\mathcal L = {\rm div} \, (A \, \nabla\cdot) +
\vec b \cdot \nabla + q$, where  
$A \in D'(\R^n)^{n\times n}$, 
 $\vec b \in D'(\R^n)^n$ and $q \in D'(\R^n)$, $n \ge 2$.  Then 
the following statements hold. 

{\rm(i)} The sesquilinear form of $\mathcal L$ is bounded, i.e., 
 {\rm(\ref{E:1.1a})} holds 
if and only  if  $A^s \in L^\infty(\R^n)^{n \times n}$, and 
 $\, \vec b$ and $q$ can be represented respectively in the form 
\begin{equation}\label{E:1.4u}
\vec b = \vec c + {\rm Div} \, F, \qquad q = {\rm div} \, \vec h,
\end{equation} 
where $F$ 
 is a skew-symmetric matrix field such that
\begin{equation}\label{E:1.4v}
F - A^c \in {\rm BMO(\R^n)}^{n\times n},
\end{equation} 
whereas $\vec c$ and $\vec h$ belong to $L^2_{{\rm loc}}(\R^n)^{n}$, and 
obey the condition 
\begin{equation}\label{E:1.5u}
|\vec c|^2 + | \vec h |^2  \in \mathfrak{M}^{1, \, 2}_+(\R^n). 
\end{equation}

{\rm(ii)} If  the sesquilinear form of $\mathcal L$ is bounded, then $\vec c$,  
$F$, and $\vec h$  in decomposition {\rm (\ref{E:1.4u})} can be determined  
explicitly by    
\begin{align}\label{E:1.6u}
\vec c & = \nabla ( \Delta^{-1} {\rm div} \, \vec b), 
 \qquad \vec h = 
 \nabla (\Delta^{-1} \, q), \\
\qquad F & = \Delta^{-1} {\rm Curl} \, [ \vec b -   {\rm Div} \, (A^c)] + 
A^c, 
\label{E:1.6uu}
\end{align} 
where 
\begin{equation}\label{E:1.6v}
\Delta^{-1} {\rm Curl} \, [ \vec b - 
{\rm Div} \, (A^c)]\in {\rm BMO}(\R^n)^{n \times n}, 
\end{equation} 
and 
\begin{equation}\label{E:1.5'u}
 \vert \nabla ( \Delta^{-1} {\rm div} \, \vec b) 
 \vert^2 + 
 \vert \nabla (\Delta^{-1} \, q)  \vert^2 \in \mathfrak{M}^{1, \, 2}_+(\R^n). 
\end{equation}} 
\end{main4}
\smallskip

  We remark that condition {\rm (\ref{E:1.6v})} in statement (ii) of Theorem IV   
may be replaced with   
\begin{equation}\label{E:1.6w}
\vec b - {\rm Div} \,  (A^c) \in {\rm BMO}^{-1}(\R^n)^{n},  
\end{equation}
which ensures that decomposition {\rm (\ref{E:1.4u})} holds. 
 Here ${\rm BMO}^{-1}(\R^n)$ stands for the space of distributions that can be represented in the form  
  $f={\rm div} \, \vec g$ where $\vec g \in {\rm BMO}(\R^n)^n$ (see \cite{KT}).

In the special case $n=2$,  it is easy to see that  {\rm(\ref{E:1.1a})} holds if 
and only if $A^s \in L^\infty(\R^2)^{2\times 2}$,  
$\vec b -  \, {\rm Div} \, (A^c) \in {\rm BMO}^{-1} (\R^2)^{2}$, and $q={\rm div} \, \vec b =0$.

As mentioned in the Introduction, expressions $\nabla (\Delta^{-1} \, q)$, 
$\nabla ( \Delta^{-1} {\rm div} \, \vec b)$, ${\rm Div} 
( \Delta^{-1} {\rm Curl} \, \vec b)$, which involve nonlocal operators,    
are defined in the sense of distributions. This is possible, since $\Delta^{-1} q$, 
 $\Delta^{-1} {\rm div} \, \vec b$, and 
$\Delta^{-1} {\rm Curl} \, \vec b$  
can be understood in terms of the convergence in the weak-$*$ topology 
of ${\rm BMO}(\R^n)$    of $\Delta^{-1} \, {\rm div} \, (\psi_N \, \vec b)$, 
$\Delta^{-1} \, {\rm Curl} \, (\psi_N \, \vec b)$, and 
$\Delta^{-1} \,  (\psi_N \, q)$, respectively, as $N \to +\infty$. 
Here  $\psi_N(x) = \psi(\frac {x} N)$ is a smooth 
cut-off function,  where $\psi$ is supported  in the unit ball $\{x: \,  |x|<1\}$, 
and $\psi (x) = 1$  if
$|x|\le \frac 1 2$. 
The limits 
above do not depend on the choice of $\psi$.

It follows from Theorem IV that $\mathcal L$  is form bounded on $L^{1, \, 2}(\R^n)\times L^{1, \, 2}(\R^n)$ if and 
only if 
$A ^s \in L^\infty(\R^n)^{n \times n}$, and $\vec b_1 \cdot \nabla +q$ 
is form bounded, where 
  \begin{equation}\label{E:1.1d}
\vec {b_1} = \vec b  - {\rm Div} (A^c). 
\end{equation}

In particular,  the principal part 
$\mathcal P u={\rm div} (A \, \nabla u)$ 
is form bounded 
if and only if  
\begin{align}\label{E:1.1e}
& A^s \in L^\infty(\R^n)^{n \times n}, \\
&  {\rm Div} \, (A^c)   
\in {\rm BMO}^{-1}(\R^n)^{n}.
\label{E:1.1ee}
\end{align} 
 A simpler condition with  $A^c \in {\rm BMO}(\R^n)^{n \times n}$ 
in place of (\ref{E:1.1ee}) is sufficient, but generally is necessary only if  $n=1,  2$. 

Thus,  the form boundedness problem for the general second 
order differential operator  
is reduced  to the special case 
\begin{equation} \label{E:1.l}
\mathcal L = \vec b\cdot \nabla  + q, \qquad \vec b \in D'(\R^n)^n, \quad q \in D'(\R^n). 
\end{equation}

  As a  corollary of Theorem~IV, 
 we deduce that, if $\vec b \cdot \nabla + q$  is form bounded, then the Hodge decomposition 
\begin{equation}\label{E:1.4uu}
\vec b = \nabla ( \Delta^{-1} {\rm div} \, \vec b) + {\rm Div} \, (\Delta^{-1} {\rm Curl} \, \vec b)  
\end{equation} 
holds, where $\Delta^{-1} ({\rm Curl} \, \vec b) \in {\rm BMO}(\R^n)^{n \times n}$, 
and 
\begin{equation}\label{E:1.7u} 
\int_{|x-y|<r}   \,  [ \,    \vert  
\nabla \Delta^{-1} ({\rm div} \, \vec b)  
\vert^2 +  \vert \nabla (\Delta^{-1} \, q)  \vert^2 \, ] \,   
dy \le {\rm const} \, \,  r^{n-2},  
\qquad 
\end{equation}
for  all 
 $r>0, \, x \in \R^n$, in the case $n\ge 3$; in two dimensions, it follows that ${\rm div}\, \vec b = q =0$.  

We observe that   condition (\ref{E:1.7u})  is generally stronger than 
 $\Delta^{-1} {\rm div} \, \vec b \in 
{\rm BMO}(\R^n)$ and $\Delta^{-1} \, q \in {\rm BMO}(\R^n)$, while  the divergence-free 
part of $\vec b$ is characterized by   
$\Delta^{-1} {\rm Curl} \, \vec b \in {\rm BMO}(\R^n)^{n \times n}$, for all $n \ge 2$.

The main difficulty in the proof of Theorem IV  is the interaction between 
the quadratic forms associated with  $q - \tfrac 1 2 \, {\rm div} \, \vec b $ and 
the divergence free part of $\vec b$. To this effect, we use Theorem II, 
which characterizes
  vector fields $\vec d$ such that the commutator 
inequality \eqref{comm-2} holds. Theorem II is proved in \cite[Lemma 4.8]{MV3} 
using the idea of 
the gauge transformation  (\cite[Sec. 7.19]{LL}; \cite[Sec. X.4]{RS}): 
$$\nabla  \to   e^{- i \lambda} \, \nabla \, e^{+ i \lambda},
$$
where the gauge $\lambda$ is 
a real-valued  function in $L^{1, \, 2}_{{\rm loc}}(\R^n)$. 

The nontrivial problem of choosing an appropriate gauge 
is solved in \cite{MV3} as follows: 
$$\lambda = \tau \, \log \, (N \mu), \qquad 1< 2 \tau < 
\tfrac n {n-2}, $$
where  $N \mu = (-\Delta)^{-1} \mu$ 
is the Newtonian potential of the equilibrium measure $\mu$
associated with an arbitrary compact set $e$ of positive  capacity. 

With this choice of $\lambda$,  the energy space $L^{1, \, 2} (\R^n)$ is  gauge 
invariant, and for the irrotational part 
  $\vec c= \nabla ( \Delta^{-1} {\rm div} \, \vec d)$  we have $|\vec c|^2  \in \mathfrak{M}^{1, \, 2}_+(\R^n)$. 
In addition, we have 
 $F=\Delta^{-1} {\rm Curl} \, \vec d$  
belongs to ${\rm BMO}(\R^n)^{n\times n}$, and $\vec d = \vec c + {\rm Div} \, F$.
  These conditions are necessary and sufficient for \eqref{comm-2}.

Applications  of Theorem IV to the 
magnetic Schr\"odinger operator $\mathcal{M}$  
defined by (\ref{E:1.1m}) are given in \cite[Theorem 3.4]{MV3}, where it is shown that  $\mathcal{M}$ is 
form bounded if and only if  both 
$q + | \vec a|^2$ and $\vec a \cdot \nabla$ are form bounded. 

\subsection{Form boundedness in $W^{1, \, 2}(\R^n) $}\label{inhomogeneous}

The above  results are easily extended to 
the Sobolev space 
$W^{1, \, 2}(\R^n)$ $(n \ge 1$) 
 with norm 
 $|| u||_{W^{1, \, 2}(\R^n)} = ||\nabla  u||_{L^2(\R^n)} + ||u||_{L^2(\R^n)}$.

In particular, necessary and sufficient 
conditions  are given in \cite[Theorem 5.1]{MV3}  
for the boundedness of the  general second order operator 
$$\mathcal{L}\!:  W^{1, \, 2}(\R^n) \to  W^{-1, \, 2}(\R^n).$$
 This solves the {\it relative form boundedness} problem  for $\mathcal{L}$,  and consequently for 
the magnetic Schr\"odinger operator $\mathcal{M}$,  
 with respect to the Laplacian  on $L^2(\R^n)$ (see \cite[Sec. X.2]{RS}). The proofs make use  of an 
 inhomogeneous version of the ${\rm div}$-${\rm curl}$ lemma 
(\cite[Lemma 5.2]{MV3}). 

\subsection{Infinitesimal form boundedness}\label{infinitesimal} 

Other fundamental properties of quadratic  
  forms associated with  differential operators  can be characterized using our methods.  
 In particular, for the Schr\"odinger operator 
$\mathcal{H} = \Delta + q$ with $q \in D'(\R^n)$, criteria of relative  compactness 
were obtained in \cite{MV1}, whereas the infinitesimal form boundedness expressed by the inequality \begin{equation}\label{E:1.9} 
| \langle q \, u, \, u \rangle | \le 
\epsilon \, ||\nabla u||^2_{L^2(\R^n)} + 
C(\epsilon) \, || u||^2_{L^2(\R^n)}, \quad   u \in C^\infty_0 
(\R^n), 
\end{equation} 
for every $\epsilon \in (0, 1)$, where $C(\epsilon)$ is a positive constant, 
along with  Trudinger's subordination where $C(\epsilon) = C \, \epsilon^{-\beta}$ $(\beta>0)$,  was characterized in \cite{MV5}. 
 Necessary and sufficient conditions for such properties in the case of the general second order differential operator  
are discussed in  \cite{MV3}.

\subsection{Nash's inequality and $p$-subordination}\label{nash} 

For $q \in D'(\R^n)$, we consider the $p$-subordination property
 \begin{equation}\label{E:1.10}
\left \vert \langle q \, u, \, u\rangle \right \vert  \leq C \, 
||\nabla u||^{2p}_{L^2(\R^n)}  \,||u||^{2(1-p)}_{L^2(\R^n)},  
\end{equation} 
for all $u 
\in C^\infty_0(\R^n)$, 
where $p  \in (0, \, 1)$.  

Nash's type inequality is similar to \eqref{E:1.10}, 
with $||u||_{L^1(\R^n)}$ in place of $||u||_{L^2(\R^n)}$ on the right-hand side,
\begin{equation}\label{E:1.11}
\left \vert \langle q \, u, \, u\rangle \right \vert  \leq C \, 
||\nabla u||^{2p}_{L^2(\R^n)}  \,||u||^{2(1-p)}_{L^1(\R^n)}.  
\end{equation} 
The classical Nash's inequality corresponds to $q\equiv 1$ 
and $p= \frac n {n+2}$ (see \cite[Theorem 8.13]{LL}. 

It is proved in \cite[Theorem 6.5]{MV5} that \eqref{E:1.10} holds  if and only if  
$
q = {\rm div} \,\,  \vec \Gamma,
$
 where $\vec \Gamma =  \nabla \Delta^{-1} q$, and one of the 
 following conditions hold:
\begin{align*}
\quad \vec \Gamma \in {\rm BMO}\quad & \hbox{if} \quad  p= 1/2;\\ 
\quad \vec \Gamma \in
{\rm Lip}({1-2p}) \quad & \hbox{if} \quad  0<p< 1/2; \\ 
 \int_{|x-y|<r} |\vec \Gamma (y)|^2 \,  dy \le c \, r^{n+2-4p}
\quad    & \hbox{if} \quad  1/2<p< 1,  
\end{align*}
 for all $r>0$ and $x \in \R^n$. Similar results hold for Nash's inequality \eqref{E:1.11} (see \cite[Corollary 6.8]{MV5}).

\subsection{Form boundedness in $W^{{\frac{1}{2}}, \, 2}(\R^n)$}\label{infinites} 

Similar problems were solved for the  fractional (modified relativistic) Schr\"odinger operator 
$\mathcal{L} = -(-\Delta)^{{\frac{1}{2}}} + q$. In particular, the boundedness of the operator 
$$\mathcal{L}\!:  W^{{\frac{1}{2}, \, 2}}(\R^n) \to  W^{-{\frac{1}{2}, \, 2}}(\R^n)$$ 
has been characterized  in \cite{MV2} using certain extensions to higher dimensions for multipliers acting from $W^{1, \, 2}(\R^{n+1})$ to $W^{-1, \, 2}(\R^{n+1})$.

\section{Accretivity}\label{accretivity}

We now turn to the accretivity problem for $-\mathcal{L}$, where $\mathcal{L}$ is a  second order linear differential operator with complex-valued distributional 
coefficients defined by \eqref{E:0.b} in an open set  $\Omega \subseteq \R^n$ ($n \ge 1$). 

\subsection{General accretivity criterion}\label{sec2.1}

Given $A=(a_{jk})\in D'(\Omega)^{n \times n}$, we define its symmetric part $A^s$ and  skew-symmetric part $A^c$  respectively by \eqref{trans}. 
The accretivity property for $-\mathcal{L}$ can be characterized  in terms of the following real-valued expressions: 
\begin{equation}\label{expr}
P={\rm Re} \, A^s, \quad  \vec{d}=  \tfrac{1}{2} \,   [ {\rm Im} \, \vec{b} - {\rm Div} \,  ( {\rm Im} \, A^c) ], \quad \sigma={\rm Re} \, c - \tfrac {1}{2}{\rm div} \, ({\rm Re} \, \vec{b}),
\end{equation}
where $P=(p_{jk})\in D'(\Omega)^{n\times n}$, $\vec{d}=({d_j})\in D'(\Omega)^{n}$, 
and $\sigma\in D'(\Omega)$. This is a consequence of the relation 
(see \cite[Sec.4]{MV6})
\begin{equation}\label{L-2}
{\rm Re}  \langle -\mathcal{L} u, u\rangle 
={\rm Re}  \langle -\mathcal{L}_2 u, u\rangle, \quad 
u\in C^\infty_0 (\Omega),
\end{equation}
where 
 \begin{equation}\label{L-2a}
  \mathcal{L}_2  = {\rm div} \, (P \nabla \cdot) 
  + 2 i \, \vec{d}\cdot \nabla + \sigma. 
 \end{equation}

Moreover, in order that $-\mathcal{L}$ be accretive,  the matrix $P$ must be nonnegative definite, i.e., $P \xi\cdot \xi \ge 0$ in 
$D'(\Omega)$ 
for all $\xi \in \R^n$. In particular, each  $p_{jj}$ ($j=1, \ldots, n$) is  a nonnegative Radon measure. 
  
  A  characterization 
  of accretive operators $-\mathcal{L}$ is given in the following criterion 
  obtained in \cite[Proposition 2.1]{MV6}. 
  
\begin{proposition}\label{prop_1}  Let $\mathcal L 
= {\rm div}  ( A \nabla \cdot) +
\vec{b} \cdot \nabla + c$, where  $A \in D'(\Omega)^{n\times n}$, 
 $\vec{b} \in D'(\Omega)^n$ and $c \in D'(\Omega)$ are complex-valued.  
 Suppose that $P$, $\vec{d}$, and $\sigma$ are defined by \eqref{expr}. 
 
 The operator $-\mathcal L$ is accretive if and only  if  $P$ is a nonnegative definite matrix, and 
the following two conditions hold: 
\begin{equation}\label{E:1.b}
[h]_{\mathcal{H}}^2=  \langle P \nabla h, \nabla h\rangle  - \langle \sigma \,  h, h\rangle\ge 0,  
\end{equation}
for all real-valued  $h \in C^\infty_0(\Omega)$,  and 
\begin{equation}\label{E:1.c}
 \left \vert \langle \vec{d},  u \nabla v -v \nabla u\rangle \right\vert 
\le [u]_{\mathcal{H}} \, [v]_{\mathcal{H}},  
\end{equation}
for all real-valued  $u, v \in C^\infty_0(\Omega)$.  
\end{proposition}

\subsection{Real-valued coefficients}\label{sec2.2} 

It follows from Proposition \ref{prop_1} that, for operators with real-valued coefficients, condition \eqref{E:1.b} alone characterizes nonnegative definite 
operators $-\mathcal{L}$ in an open set $\Omega\subseteq \R^n$ ($n\ge 1$). 
   A more explicit characterization of  this property,  
  under 
 the assumption that  $P=A^s\in L^1_{{\rm loc}} (\Omega)^{n \times n}$  in the sufficiency part, and 
 that $P$ is uniformly elliptic in the necessity part, is given in the next proposition 
 (see \cite[Theorem 2.2]{MV6}). 
 
\begin{proposition}\label{prop_2}
Let $\mathcal{L} = {\rm div}  ( A \nabla \cdot) +
\vec{b} \cdot \nabla + c$, where  $A \in D'(\Omega)^{n\times n}$, 
 $\vec{b} \in D'(\Omega)^n$ and $c \in D'(\Omega)$ are real-valued. 
 Suppose that $P=A^s\in L^1_{{\rm loc}} (\Omega)^{n \times n}$ is a nonnegative definite matrix a.e. 
 
 {\rm (i)} If  there exists a measurable vector field $\vec{g}$ in $\Omega$ such that  $(P \vec{g}) \cdot \vec{g} \in L^1_{{\rm loc}}(\Omega)$, and 
 \begin{equation}\label{g-cond}
   \sigma=c - \tfrac {1}{2}{\rm div} \, (\vec{b}) \le {\rm div} \, (P \vec{g}) 
   - (P \vec{g}) \cdot \vec{g} \quad \textrm{in} \, \, D'(\Omega), 
      \end{equation}
then the operator 
$-\mathcal{L}$ is nonnegative definite. 

 {\rm (ii)} Conversely, if $-\mathcal{L}$ is nonnegative definite,  then there exists a vector field $\vec{g} \in L^2_{{\rm loc}} (\Omega)^n$ so that 
 $(P \vec{g}) \cdot \vec{g} \in L^1_{{\rm loc}}(\Omega)$, and 
 \eqref{g-cond} holds, 
    provided $P$ is uniformly elliptic.  
\end{proposition}

The uniform ellipticity condition on $P$ in statement (ii) of Proposition \ref{prop_2} can be relaxed. This question will be treated elsewhere. 

Results  similar to Proposition \ref{prop_2}  are well known
 in ordinary differential equations  \cite[Sec. XI.7]{Ha}, in relation  to disconjugate Sturm-Liouville equations and Riccati equations with continuous coefficients  
 (see also
 \cite{H}, \cite{MV2}, \cite{MV6}).

\subsection{Nonnegative definite Schr\"{o}dinger operators}\label{sec2.3} 

As was mentioned above, in the special case  of Schr\"{o}dinger operators 
$ \mathcal{H}={\rm div} \, (P \nabla h) + \sigma$,  with real-valued  $\sigma \in D'(\Omega)$ and uniformly elliptic 
$P$, Proposition \ref{prop_2} was obtained originally in 
\cite[Proposition 5.1]{JMV1}. Under these assumptions, $-\mathcal{H}$ is nonnegative definite, i.e., 
\[
[h]_{\mathcal{H}}^2=\langle -\mathcal{H} h, h \rangle\ge 0, \quad \textrm{for all} \, \, h \in C^\infty_0(\Omega), 
\]
 if and only if there exists a vector field $\vec{g} \in 
        L^2_{{\rm loc}} (\Omega)^n$ such that 
    \begin{equation}\label{g-cond1}
   \sigma \le {\rm div} \, (P \vec{g}) 
   - P \vec{g} \cdot \vec{g} \quad \textrm{in} \, \, D'(\Omega).
      \end{equation}

A simpler linear sufficient condition for $-\mathcal{H}$ to be nonnegative definite is given by 
$
 \sigma \le {\rm div} \, (P \vec{g}), 
$
where $\vec{g} \in L^2_{{\rm loc}} (\Omega)^n$ satisfies the inequality
\[
\int_\Omega (P \vec{g} \cdot \vec{g}) \, h^2 \, dx \le \frac{1}{4}\int_\Omega 
|P \nabla h|^2 \, dx, \quad \textrm{for all} \, \, h \in C^\infty_0(\Omega).
\]
 Here $P \vec{g} \cdot \vec{g} \in  \mathfrak{M}^{1, \, 2}_+(\Omega)$, and so 
 $|\vec{g}|^2$ is admissible  if $P$ is uniformly elliptic. However, 
such conditions are not necessary, with any constant in place of $ \frac{1}{4}$, even when $P=I$; see 
\cite[Proposition 7.1]{JMV1}.

We observe that in Proposition \ref{prop_1} above, 
 the  nonnegative definite quadratic form $[h]_{\mathcal{H}}^2$ is associated with  the Schr\"{o}dinger operator $-\mathcal{H}$, where  $\mathcal{H}$ has real-valued coefficients  
 $
  P ={\rm Re} \, A^s$ and  
$ \sigma={\rm Re} \, c - \frac {1}{2}{\rm div} \, ({\rm Re} \, \vec{b}) $. 
 Hence,  \eqref{g-cond1} characterizes 
      the first condition of Proposition \ref{prop_1}  given by 
      \eqref{E:1.b}. The second one, namely, the commutator condition \eqref{E:1.c}, will be discussed further in Sections \ref{sec2.5} 
      and \ref{sec2.6}.

\subsection{The one-dimensional case}\label{sec2.4} 

In this section, the 
differential operator  $\mathcal{L} u = (a \, u')' + b u' + c$ is  defined on an open interval $I \subseteq \R$ (possibly unbounded).
 In this case, one can avoid commutator estimates 
using methods  of ordinary differential equations (\cite{Ha}, \cite{Hi}).  In the statements below we will make use of the standard convention $\frac{0}{0}=0$.  The following criterion of accretivity for complex-valued coefficients in the one-dimensional case was obtained in \cite[Theorem 2.2]{MV6}.

\begin{main5}\label{Theorem V} {\it  Let $a, b, c  \in D'(I)$. Suppose that $p={\rm Re} \, a  \in L^1_{{\rm loc}} (I)$, and ${\rm Im} \,  b \in L^1_{{\rm loc}} (I)$. 

{\rm (i) } The operator $-\mathcal{L}$ is accretive  
 if and only if  $\frac{({\rm Im} \,  b)^2}{p}\in L^1_{{\rm loc}}(I)$, where $p \ge 0$ a.e.,
and the following quadratic form inequality holds:
\begin{equation}\label{0.1}
\int_{I} p (h')^2 dx -\langle {\rm Re} \, c - \frac{1}{2}({\rm Re} \, b)',  h^2 
\rangle - \int_{I} \frac{({\rm Im} \, b)^2}{4p}\, h^2 \, dx \ge 0,
\end{equation}
for all real-valued $h \in C^{\infty}_0(I)$.

{\rm (ii) }    If  there exists a function 
$f \in L^1_{{\rm loc}}(I)$ such that $ \frac{f^2}{p} \in L^1_{{\rm loc}}(I)$, and 
\begin{equation}\label{0.2}
{\rm Re} \, c - \frac{1}{2}({\rm Re} \, b)' 
- \frac{({\rm Im} \, b)^2}{4 p}\le f' - \frac{f^2}{p} \quad \textrm{in} \, \, D'(I), 
\end{equation}
then the operator $-\mathcal{L}$ is accretive. 

Conversely, if $-\mathcal{L}$ is accretive, and 
$m\le p(x)\le M$ a.e. for some constants $M, m>0$, then 
there exists a function 
$f \in L^2_{{\rm loc}}(I)$ such that \eqref{0.2} holds.} 
\end{main5}
\smallskip

We remark that in Theorem V, the terms 
${\rm Im} \, a$ and ${\rm Im} \, c$ play no role, but the behavior of 
${\rm Im} \, b$ is essential.  In higher dimensions, the situation is even more complicated.  
 The term ${\rm Im} \, b$ may contain both the irrotational and divergence-free 
 components, and the latter may interact 
with  
${\rm Im} \, A^c$.

\subsection{Upper  and lower bounds of quadratic forms}\label{sec2.5}

For general operators with complex-valued coefficients in the case $n \ge 2$, we recall that the first condition of Proposition \ref{prop_1}  is  necessary 
 for the accretivity of $-\mathcal{L}$, namely, 
\begin{equation}\label{nec}
  \langle  \sigma \,  h, h\rangle 
\le  \int_\Omega (P \nabla h \cdot \nabla h) \, dx,  
\end{equation}
for all real-valued $h \in C^\infty_0(\Omega)$, where $\sigma = {\rm Re } \, c - \frac{1}{2}{\rm div} ({\rm Re} \,  \vec{b}) \in D'(\Omega)$, and 
${\rm Re} \, A^{s}=P \in D'(\Omega)^{n \times n} $ is a nonnegative definite matrix. 

Suppose now that $\sigma$
has a slightly smaller upper form bound, that is, 
\begin{equation}\label{E:ia}
  \langle  \sigma \,  h, h\rangle 
\le (1-\epsilon^2)  \int_\Omega (P \nabla h\cdot \nabla h) \, dx, \quad h \in C^\infty_0(\Omega),  
\end{equation}
for some $\epsilon \in (0, 1]$. We also consider the corresponding lower bound, 
\begin{equation}\label{E:ib}
  \langle  \sigma \,  h, h\rangle 
\ge - K  \int_\Omega (P \nabla h\cdot \nabla h) \, dx, \quad h \in C^\infty_0(\Omega),  
\end{equation}
for some  constant $K\ge 0$. 

Such restrictions on  real-valued $\sigma \in D'(\Omega)$ were 
invoked in \cite[Theorem 1.1]{JMV1}, for uniformly elliptic $P$.  

We observe that 
\eqref{E:ia} is  satisfied for any $\epsilon \in (0, 1)$, 
up to an extra term  $C \, ||h||^2_{L^2(\Omega)}$, 
 if 
$\sigma$ is \textit{infinitesimally form bounded}  (see Sec. \ref{infinitesimal}).   The second term  on the right is sometimes  included in the definition of accretivity of the operator $-\mathcal{L}$. We can always incorporate it 
as a constant term in $\sigma -C(\epsilon)$. The same is true with regards to the lower bound where we can use $\sigma +C(\epsilon)$.

Assuming that both bounds \eqref{E:ia} and \eqref{E:ib} hold for some $\epsilon\in (0, 1]$ and $K\ge 0$, we obviously have, for all  $h \in C^\infty_0(\Omega)$, 
\begin{equation}\label{K-epsilon}
\epsilon   \int_\Omega (P \nabla h\cdot \nabla h) \, dx \le [h]^2_{\mathcal{H}} \le (K+1)^{\frac{1}{2}}  \int_\Omega ( P \nabla h \cdot \nabla h) \, dx. 
\end{equation}

If $P$ satisfies the uniform ellipticity assumptions \eqref{ell}, then from 
\eqref{K-epsilon} it follows that 
condition \eqref{E:1.c} equivalent, up to a constant multiple, to 
\begin{equation}\label{E:i.e}
\left \vert \langle \vec{d},  u \nabla v -v \nabla  u\rangle \right\vert 
\le C \,  ||\nabla u||_{L^2(\Omega)} \, ||\nabla v||_{L^2(\Omega)}
\end{equation}
where $C>0$ is a constant which does not depend on real-valued $u, v \in C^\infty_0(\Omega)$. For $\Omega=\R^n$ and $\vec{d}\in L^1_{{\rm loc}}(\R^n)$, see Theorem II above. 

In the case $\Omega=\R^n$, inequality \eqref{E:i.e} was characterized completely 
 in \cite[Lemma 4.8]{MV3} 
for complex-valued $u, v$. However, that characterization obviously works  
in the case of real-valued $u, v$ as well (one only needs to change the constant $C$ up to a factor of $\sqrt{2}$).

\subsection{Accretivity criterion in $\R^n$}\label{sec2.6} 

Combining the characterization of the commutator inequality \eqref{E:i.e} 
 with 
Proposition \ref{prop_1} yields the following accretivity criterion 
(\cite[Theorem 2.7]{MV6}), where the lower  bound \eqref{E:ib}  in used the necessity part, whereas the upper bound \eqref{E:ia}  
is invoked in the sufficiency part.

\begin{main6}\label{Theorem VI} {\it Let  $\mathcal L$ be the second order differential 
operator  \eqref{E:0.b} on $\R^n$ $(n \ge 2)$ with complex-valued coefficients $A \in D'(\R^n)^{n\times n}$, 
 $\vec{b} \in D'(\R^n)^n$ and $c \in D'(\R^n)$.  Let $P$, 
 $\vec{d}$ and $\sigma$ be defined by \eqref{expr}, where $P$ is uniformly 
 elliptic. 
 
  {\rm (i) } Suppose that $-\mathcal{L}$ is accretive, i.e.,  \eqref{accr}  holds, and 
  $\sigma$ satisfies \eqref{E:ib}  for some $K\ge 0$. Then  $\vec{d}$  can be represented in the form 
\begin{equation}\label{E:ig}
\vec{d} = \nabla f + {\rm Div} \, G, 
\end{equation} 
where $f\in D'(\R^n)$ is real-valued,  $|\nabla f |^2 \in \mathfrak{M}^{1, \, 2}_+(\R^n)$,  
and $G \in {\rm BMO}(\R^n)^{n\times n}$ is 
a real-valued skew-symmetric matrix field. 

Moreover,   
$f$ and  $G$ above can be defined explicitly as 
\begin{equation}\label{E:igb}
f = \Delta^{-1} ( {\rm div } \, \vec{d}), \qquad 
G =\Delta^{-1} ( {\rm Curl} \, \vec{d}).
\end{equation}

{\rm (ii) } Conversely, suppose that  $\sigma$ satisfies \eqref{E:ia} 
with some $\epsilon \in (0, 1]$. Then $-\mathcal{L}$ is accretive  
if representation \eqref{E:ig} holds, where $|\nabla f |^2 \in \mathfrak{M}^{1, \, 2}_+(\R^n)$,  
and $G \in {\rm BMO}(\R^n)^{n\times n}$ is 
a real-valued skew-symmetric matrix field, provided both 
$\Vert |\nabla f|^ 2 \Vert_{\mathfrak{M}^{1, \, 2}_+(\R^n)}$ and the 
${\rm BMO}$-norm of $G$ 
are small enough,  
depending only on  $\epsilon$.} 
\end{main6}
\smallskip

If $n=2$, then  in Theorem VI, we have $f=0$,  and $\vec{d}= (-\partial_2 g, \partial_1 g)$
with $g \in {\rm BMO}(\R^2)$. In statement (ii), the ${\rm BMO}$-norm of $g$ is 
supposed to be  small enough
(depending only on  $\epsilon$).

  If $n=3$, one can use the usual vector-valued ${\rm curl} (\vec{g}) \in D'(\R^3)^3$ in place of ${\rm Div} \, G$ in decomposition \eqref{E:ig}, with   
$\vec{g} =\Delta^{-1} ( {\rm curl} \, \vec{d})$ in place of $G$  in \eqref{E:igb}.


\end{document}